
  \advance \hsize by 2pt

  \newcount\fontset
  \fontset=1
  \def\dualfont#1#2#3{\font#1=\ifnum\fontset=1 #2\else#3\fi}

  \dualfont\bbfive{bbm5}{cmbx5}
  \dualfont\bbseven{bbm7}{cmbx7}
  \dualfont\bbten{bbm10}{cmbx10}
  \font \eightbf = cmbx8
  \font \eighti = cmmi8 \skewchar \eighti = '177
  \font \eightit = cmti8
  \font \eightrm = cmr8
  \font \eightsl = cmsl8
  \font \eightsy = cmsy8 \skewchar \eightsy = '60
  \font \eighttt = cmtt8 \hyphenchar\eighttt = -1
  
  \font \sixbf = cmbx6
  \font \sixi = cmmi6 \skewchar \sixi = '177
  \font \sixrm = cmr6
  \font \sixsy = cmsy6 \skewchar \sixsy = '60
  \font \tensc = cmcsc10
  
  \font \titlefont = cmr7 scaled \magstep4
  \scriptfont \bffam = \bbseven
  \scriptscriptfont \bffam = \bbfive
  \textfont \bffam = \bbten

  \font\rs=rsfs10 

  \newskip \ttglue

  \def \eightpoint {\def \rm {\fam0 \eightrm }%
  \textfont0 = \eightrm
  \scriptfont0 = \sixrm \scriptscriptfont0 = \fiverm
  \textfont1 = \eighti
  \scriptfont1 = \sixi \scriptscriptfont1 = \fivei
  \textfont2 = \eightsy
  \scriptfont2 = \sixsy \scriptscriptfont2 = \fivesy
  \textfont3 = \tenex
  \scriptfont3 = \tenex \scriptscriptfont3 = \tenex
  \def \it {\fam \itfam \eightit }%
  \textfont \itfam = \eightit
  \def \sl {\fam \slfam \eightsl }%
  \textfont \slfam = \eightsl
  \def \bf {\fam \bffam \eightbf }%
  \textfont \bffam = \eightbf
  \scriptfont \bffam = \sixbf
  \scriptscriptfont \bffam = \fivebf
  \def \tt {\fam \ttfam \eighttt }%
  \textfont \ttfam = \eighttt
  \tt \ttglue = .5em plus.25em minus.15em
  \normalbaselineskip = 9pt
  \def \MF {{\manual opqr}\-{\manual stuq}}%
  \let \sc = \sixrm
  \let \big = \eightbig
  \setbox \strutbox = \hbox {\vrule height7pt depth2pt width0pt}%
  \normalbaselines \rm }


  \def \Headlines #1#2{\nopagenumbers
    \voffset = 2\baselineskip
    \advance \vsize by -\voffset
    \headline {\ifnum \pageno = 1 \hfil
    \else \ifodd \pageno \tensc \hfil \lcase {#1} \hfil \folio
    \else \tensc \folio \hfil \lcase {#2} \hfil
    \fi \fi }}

  \def \Title #1{\vbox{\baselineskip 20pt \titlefont \noindent #1}}

  \def \Date #1 {\footnote {}{\eightit Date: #1.}}

  \def \Authors #1{\bigskip \bigskip \noindent #1}

  \long \def \Addresses #1{\begingroup \eightpoint \parindent0pt
\medskip #1\par \par \endgroup }

  \long \def \Abstract #1{\begingroup \eightpoint
  \bigskip \bigskip \noindent
  {\sc ABSTRACT.} #1\par \par \endgroup }


  \def \lcase #1{\edef \auxvar {\lowercase {#1}}\auxvar }
  \def \vg #1{\ifx #1\null \null \else
    \ifx #1\ { }\else
    \ifx #1,,\else
    \ifx #1..\else
    \ifx #1;;\else
    \ifx #1::\else
    \ifx #1''\else
    \ifx #1--\else
    \ifx #1))\else
    { }#1\fi \fi \fi \fi \fi \fi \fi \fi \fi }

  \def \goodbreak {\vskip0pt plus.1\vsize \penalty -250 \vskip0pt
plus-.1\vsize }

  \newcount \secno \secno = 0
  \newcount \stno

  \def \seqnumbering {\global \advance \stno by 1
    \number \secno .\number \stno }

  \def \label #1{\def\localvariable {\number \secno
    \ifnum \number \stno = 0\else .\number \stno \fi }\global \edef
    #1{\localvariable }}

  \def\section #1{\global\def\SectionName{#1}\stno = 0 \global
\advance \secno by 1 \bigskip \bigskip \goodbreak \noindent {\bf
\number \secno .\enspace #1.}\medskip \noindent \ignorespaces}

  \long \def \sysstate #1#2#3{\medbreak \noindent {\bf \seqnumbering
.\enspace #1.\enspace }{#2#3\vskip 0pt}\medbreak }
  \def \state #1 #2\par {\sysstate {#1}{\sl }{#2}}
  \def \definition #1\par {\sysstate {Definition}{\rm }{#1}}


  \def \proof {\medbreak \noindent {\it Proof.\enspace }}
  \def \proofend {\ifmmode \eqno \square \else \hfill \square
\looseness = -1 \medbreak \fi }

  \def \$#1{#1 $$$$ #1}
  \def\=#1{\buildrel (#1) \over =}

  \def\iItem {\smallskip}
  \def\Item #1{\smallskip \item {#1}}
  \newcount \zitemno \zitemno = 0
  \def\izitem {\zitemno = 0}
  \def\zitem {\global \advance \zitemno by 1 \Item {{\rm(\romannumeral
\zitemno)}}}

  \newcount \footno \footno = 1
  \newcount \halffootno \footno = 1
  \def\footcntr {\global \advance \footno by 1
  \halffootno =\footno
  \divide \halffootno by 2
  $^{\number\halffootno}$}


  \def \({\left (}
  \def \){\right )}
  \def \[{\left \Vert }
  \def \]{\right \Vert }
  \def \*{\otimes }
  \def \+{\oplus }
  \def \:{\colon }
  \def \<{\left \langle }
  \def \>{\right \rangle }
  \def \text #1{\hbox {\rm #1}}
  \def \and {\hbox {,\quad and \quad }}

  \def \cstar {$C^*$}
  
  \def \inv {^{-1}}
  
  \def \square {\hbox {$\sqcap \!\!\!\!\sqcup $}}
  \def \stress #1{{\it #1}\/}
  
  \def \x {\times }
  \def \|{\Vert }
  \def \inv {^{-1}}


  \def\cite #1{{\rm [\bf #1\rm ]}}
  \def\scite #1#2{\cite{#1{\rm \hskip 0.7pt:\hskip 2pt #2}}}
  \def\lcite #1{(#1)}
  
  \def\bibitem#1#2#3#4{\smallskip \item {[#1]} #2, ``#3'', #4.}

  \def \references {
    \begingroup
    \bigskip \bigskip \goodbreak
    \eightpoint
    \centerline {\tensc References}
    \nobreak \medskip \frenchspacing }

  \font\titlefont = cmr7 scaled \magstep5
  \font\rs=rsfs10
  \font\forsmallbullet=cmsy10 scaled 500 

  \def\Lin{\hbox{\rs B}}
  \def\smallbullet{\raise 1pt \hbox{\forsmallbullet }}
  
  \def\B{{\cal B}}
  \def\J{{\cal J}}
  \def\Ker{{\rm Ker}}
  \def\({\big(}
  \def\){\big)}
  \def\Cr#1{C^*_r(#1)}
  \def\CrB{\Cr\B}
  \def\CG{\Cr G}
  
  \def\vflx#1#2{\mathrel {\buildrel \hbox{$#1$} \over {\hbox to #2em
{\rightarrowfill}}}}
  \def\flx#1{\vflx{#1}{1.5}}
  \def\flxb{\big\downarrow}
  \def\flxc{\big\uparrow}
  \def\sse{\ \Leftrightarrow\ }
  \def\Ideal{{\cal I}}

  \def\BMS{BMS}
  \def\FD{FD}
  \def\Dix{D}
  \def\Amena{E2}
  \def\uncondit{E1}
  \def\EL{EL}
  \def\Gromov{G}
  \def\Ozawa{O}
  \def\SV{SV}
  \def\Nica{N}
  \def\KW{KW}
  \def\topfree{EL}
  \def\Wasserman{W}

  \def\titletext{Exact Groups, Induced Ideals, and Fell Bundles}

  \Headlines {\titletext} {Ruy Exel}

  \Title
  {\titletext}

  \Date {August 28, 2000}

  \Authors
  {Ruy Exel\footnote{*}{\eightrm Partially supported by CNPq.}}

  \Addresses
  {Departamento de Matem\'atica\par
  Universidade Federal de Santa Catarina\par
  88040-900 Florian\'opolis SC\par
  BRAZIL\vskip 4pt
  E-mail: exel@mtm.ufsc.br}

  \Abstract {Given a C*-algebra $B$ which is graded over a discrete
group $G$ we consider ideals of $B$ which are invariant under the
projections onto each of the grading subspaces.  If $G$ is exact and
the standard conditional expectation of $B$ is faithful we show that
all such ideals are induced, i.e.~are generated by their intersection
with the unit fiber algebra.  This result is derived from a
generalization to the context of Fell bundles of a Theorem by
Kirchberg and Wasserman according to which a group is exact if and
only if its reduced C*-algebra is exact.}

  \section{Introduction}
  Recall from \cite{\KW} that a locally compact group $G$ is said to
be exact if the operation of taking reduced crossed products by $G$
preserves short exact sequences of C*-algebras.  More specifically $G$
is exact if, given any short exact sequence
  $$
  0 \to J \to A \to B \to 0
  $$
  of C*-algebras carrying actions of $G$, and such that all maps are
covariant, the sequence resulting from taking the reduced crossed
product by $G$ remains exact.  In \scite{\KW}{5.2} Kirchberg and
Wasserman showed that a discrete group $G$ is exact in the above sense
if and only if its reduced C*-algebra $C^*_r(G)$ is an exact
C*-algebra (i.e.~if taking the minimal tensor product with $C^*_r(G)$
preserves short exact sequences of C*-algebras \cite{\Wasserman}).

It has long been realized that crossed products may be viewed as a
special case
  of cross-sectional C*-algebras in the context of Fell bundles (see
\cite{\FD} for a comprehensive treatment of the theory of Fell
bundles, also referred to as \cstar-algebraic bundles) and hence it is
natural to ask for a generalization of Kirchberg and Wasserman's
result to this wider context.  It is the purpose of the present work
to prove precisely such a generalization.

As a byproduct we obtain a result about induced ideals on C*-algebras
which are graded over exact groups along the lines of similar results
obtained by Str\v{a}til\v{a} and Voiculescu on AF-algebras \cite{\SV},
by Nica on quasi-lattice ordered groups \cite{\Nica}, and by the
author on Fell bundles \cite{\Amena}.  See also Remark 3.2 in
\cite{\topfree}.

We say that a C*-algebra $B$ is graded over a discrete group $G$ if
$B$ is the closure of a direct sum of the form $\bigoplus_{t\in G}
B_t$ where the $B_t$'s are closed linear subspaces satisfying $B_tB_s
\subseteq B_{ts}$ and $B_t^*=B_{t\inv}$.  In all of the interesting
examples there is a conditional expectation of $B$ onto $B_e$
vanishing on the other $B_t$'s.  By \scite{\Amena}{3.5} the existence
of such a map implies that the canonical projections
  $$
  F_t : \bigoplus_{t\in G} B_t \to B_t
  $$
  extend to bounded linear maps on $B$.  For an element $b\in B$ the
\stress{coefficients} $F_t(b)$ often play the role of Fourier
coefficients and much of the harmonic analysis regarding the
convergence of Fourier series carry over to this situation.  It is
therefore natural to make the hypothesis that such a conditional
expectation exists and in such a case we say
(cf.~\scite{\Amena}{Definition 3.4})
that the grading is a \stress{topological grading}.

We say that an ideal $\Ideal$ in $B$ is \stress{invariant} if
$F_t(\Ideal)\subseteq \Ideal$ for all $t$ in $G$ (the reason for this
terminology is that when $G$ is abelian this is equivalent to $\Ideal$
being invariant under the dual action, that is the action of the
Pontryagin dual $\widehat G$ on $B$ given by $\alpha_x(b_t) = \<x,t\>
b_t$ for all $t\in G$, $x\in \widehat G$, and $b_t\in B_t$
\scite{\uncondit}{5.1}).

Following \scite{\Amena}{3.10} we say that $\Ideal$\/ is an
\stress{induced} ideal when $\Ideal = \<\Ideal \cap B_e\>$.  If $G$ is
an amenable group it is possible to show that every invariant ideal is
induced by showing that any element $b\in\Ideal$ can be recovered from
its Fourier coefficients, if not by summing its Fourier series, at
least by a Cesaro means argument \scite{\Amena}{4.9 and 4.10}.

On the other hand, for any non-amenable group one can produce an
example of an invariant ideal which is not induced: just take the
kernel of the regular representation of $C^*(G)$.  This seems to
indicate that one should concentrate on ``reduced'' gradings, meaning
topological gradings whose conditional expectation is faithful.
Nevertheless, given that non-exact discrete groups exist (see
\cite{\Ozawa}, \cite{\Gromov}), even under the hypothesis that the
grading is reduced, one may construct an example of an invariant ideal
which is not induced, as is done at the end of Section 4 in
\cite{\Amena}.

In the main application of the present work we deduce, as a
consequence of our generalization of Kirchberg and Wasserman's
result, that invariant
ideals are necessarily induced whenever the base group is exact and
the grading is reduced.

Last but not least we would like to acknowledge a few short but
fruitful discussions with Simon Wasserman and Iain Raeburn from where
some of the techniques used in this work evolved.

  \section{Ideals in Fell Bundles}
  Throughout this work we shall let $G$ be a discrete group and $\B =
\{B_t\}_{t\in G}$ be a Fell bundle over $G$. 

  In this section we will briefly outline the main aspects of the
theory of ideals in Fell bundles which we will need in the sequell.
We shall leave out most proofs since they are always easy
generalizations of the corresponding ones for C*-algebras.

  \definition 
  \label \DefineIdeal
  An \stress{ideal} in $\B$ is a colection $\J = \{J_t\}_{t\in G}$,
where each $J_t$ is a closed subspace of $B_t$, such that for all
$s,t\in G$ one has that
  $B_s J_t \subseteq J_{st}$ and
  $J_s B_t \subseteq J_{st}$.

For example, if $J$ is an ideal in the cross-sectional C*-algebra
$C^*(\B)$ then, defining $J_t = J\cap B_t$, one has that
$\{J_t\}_{t\in G}$ is  an ideal in $\B$.

Given an ideal $\J$ of $\B$ one clearly has that $J_e$ ($e$ denoting
the unit element in $G$) is an ideal of $B_e$.  If $\{u_i\}$ is an
approximate identity for $J_e$ it is easy to show that for each $b$ in
each $J_t$ one has that $\lim_i u_i b = \lim_i b u_i= b$.  As in the
case of C*-algebras \scite{\Dix}{1.8.2} it follows that $\J$ is
self-adjoint in the sense that $(J_t)^*=J_{t\inv}$.  One then has that
$\J$ is a Fell bundle in its own right.

Consider, for each $t$ in $G$, the quotient space $B_t/J_t$.  It is
clear that the operations of $\B$ drop to the quotient giving a
multiplication operation
  $$
  \smallbullet : {B_s \over J_s} \times {B_t \over J_t}\ \to \ {B_{st}
\over J_{st}}
  $$
  and an involution
  $$
  * : {B_t \over J_t}\ \to \ {B_{t\inv}\over J_{t\inv}}.
  $$
  With the quotient topology one has that the collection $\B/\J :=
\{B_t/J_t\}_{t\in G}$ is then a Fell bundle (the proof of the
C*-identity $\|b^*b\| = \|b\|^2$, for $b\in B_t/J_t$, is perhaps the only
slightly nontrivial verification to be made but it again follows as in
\scite{\Dix}{1.8.2}).

The inclusion map $\J\to\B$ and the quotient map $\B\to \B/\J$ clearly
give *-homomor\-phisms 
  $$
  C^*(\J)\buildrel \iota \over \to C^*(\B) {\rm\quad and\quad }
  C^*(\B)\buildrel \kappa \over \to C^*(\B/\J).
  $$

  \state Proposition 
  \label \FullExactSeq
  If $\J$ is an ideal in $\B$ then
  $$
  0 \to C^*(\J) \flx \iota C^*(\B) \flx \kappa C^*(\B/\J) \to 0
  $$
  is an exact sequence of C*-algebras.

  \proof
  Let $\pi$ be a faithful non-degenerate representation of $C^*(\J)$
on a Hilbert space $H$.  By \scite{\FD}{VIII.9.4} we have that
$\pi|_{J_e}$ is also non-degenerate.

Let $\{u_i\}$ be an approximate identity for $J_e$.  It is easy to see
that for each $b_t$ in each $B_t$ the net
  $
  \{\pi(u_i b_t)\}_i
  $
  converges in the strong operator topology.  Denoting the limit by
$\widetilde \pi(b_t)$ one clearly has that
  $\widetilde \pi(b_t) \pi(a) \xi = \pi(b_t a)\xi$
  for all $a\in J_e$ and $\xi\in H$, which in turn can be used to
prove that $\widetilde \pi$ is a *-representation of $\B$ in the sense
of \scite{\FD}{VIII.9}.  Denote also by $\widetilde \pi$ the integrated
form of $\widetilde \pi$ \scite{\FD}{VIII.11.6} which is then a
representation of $C^*(\B)$ extending $\pi$.  Given $a\in C^*(\J)$ we
therefore have that
  $$
  \|a\| = \|\pi(a)\| = \|\widetilde\pi(\iota(a))\| \leq \|\iota(a)\|
  $$
  thus proving that $\iota$ is one-to-one.

The composition $\kappa\circ\iota$ clearly vanishes and hence the range
of $\iota$ is contained in $\Ker(\kappa)$.  In order to prove the
reverse containment
observe that, for each $t\in G$, the composition
  $$
  B_t \to C^*(\B) \buildrel q \over \to C^*(\B)/C^*(\J),
  $$
  where the leftmost arrow is the canonical inclusion of $B_t$ in
$C^*(\B)$ and $q$ is the canonical quotient map, vanishes on $J_t$ and
hence drops to the quotient yielding a map
  $$
  j_t : B_t/J_t \to  C^*(\B)/C^*(\J).
  $$
  Together the collection of maps $\{j_t\}_{t\in G}$ forms a
representation of $\B/\J$ into $C^*(\B)/C^*(\J)$ which in turn gives
rise to an integrated representation, namely a *-homomorphism
  $$
  \psi : C^*(\B/\J) \to C^*(\B)/C^*(\J),
  $$
  whose restriction to $B_t/J_t$ coincides with $j_t$.  One may now
easily show that $\psi\circ\kappa$ is the canonical quotient map $q:
C^*(\B) \to C^*(\B)/C^*(\J)$.  Thus, if $a\in\Ker(\kappa)$ we have that
$0 = \psi(\kappa(a))=q(a)$ and hence $a\in C^*(\J)$ as desired.  This
shows that our sequence is exact at the middle.  Exactness at
$C^*(\B/\J)$, namely the fact that $\kappa$ is onto, is left as an easy
exercise.
  \proofend

  \section{The absorption property of the left regular representation}
As before we will let $\B = \{B_t\}_{t\in G}$ be a Fell bundle over
the discrete group $G$.  
  In this section we aim to prove that the tensor product (in a sense
to be made precise below) of any representation of $C^*(\B)$ with the
left regular representation of $G$ leads to a representation which
factors through the reduced cross-sectional C*-algebra $\CrB$.

We refer the reader to \cite{\Amena} for the basic theory of reduced
cross-sectional C*-algebras of Fell bundles but let us   nevertheless
recall a couple of facts from \cite{\Amena} which will be crucial in
what follows:

  \state Proposition
  \label \Recollection
  There exists a *-homomorphism
  $\Lambda: C^*(\B)\to \CrB$
  which is the identity on each $B_t$ and whose kernel is given  by
  $$
  \Ker(\Lambda) = \Big\{b \in C^*(\B) : E(b^*b)=0\Big\},
  $$
  where $E$ is the standard conditional expectation from $C^*(\B)$ to
$B_e$.

  \proof
  See the paragraph immediately after \scite{\Amena}{2.3} as well as
\scite{\Amena}{3.6}.
  \proofend

Let $\Pi$ be any representation of $C^*(\B)$ on a Hilbert space $H$,
considered fixed for the time being.  The restriction of $\Pi$ to the
disjoint union of the $B_t$'s is easily seen to form a representation
of $\B$
  which we shall denote by $\pi$.  Consider the representation
$\pi\*\lambda$ of $\B$ on $H\*\ell_2(G)$ given by
  $$
  (\pi\*\lambda)(b_t) = \pi(b_t)\*\lambda_t,
  $$
  whenever $t\in G$ and $b_t\in B_t$ (as usual the ``$\lambda$'' in
the right hand side above refers to the left regular representation of
$G$).

\definition We shall denote by $\Pi\*\Lambda$ the representation of
$C^*(\B)$ obtained by integrating $\pi\*\lambda$.

Denote by $\{\delta_t\}_{t\in G}$ the canonical orthonormal basis of $\ell_2(G)$
and for each $t$ let $P_t$ be the orthogonal projection from $H\*\ell_2(G)$ onto
$H\*\delta_t$.  In what follows we shall make use of the \stress{compression
operator}
  $$
  F_t\ :\ \Lin\(H\*\ell_2(G)\) \ni T\ \longmapsto\ 
  P_tTP_t \in \Lin(H\*\delta_t).
  $$

  \state Lemma
  \label \FIsFaithful
  If\/ $T$ is a nonnegative operator belonging to the range of\/
$\Pi\*\Lambda$ and $F_e(T)=0$ then $T=0$.

  \proof
  Denote by $\rho$ the right regular representation of $G$ on
$\ell_2(G)$, so that
  $\rho_t(\delta_s) = \delta_{st\inv}$  for $t,s\in G$.
  Observe that $1\*\rho_t$ commutes
with the range of $\Pi\*\Lambda$ while 
  $$(1\*\rho_{t\inv}) P_s (1\*\rho_t) = P_{st},$$
  for all $s$.  So, for any $T$ in the range of $\Pi\*\Lambda$, we obtain
  $$
  (1\*\rho_{t\inv}) F_s(T) (1\*\rho_t) = 
  F_{st}(T).
  $$
  Assuming that $F_e(T)=0$ we then conclude that  $F_t(T)=0$  for all
$t$ in $G$ which therefore says that the diagonal of $T$ vanishes 
relative to the decomposition
  $$H\*\ell_2(G) = \bigoplus_{t\in G} H\*\delta_t.$$
  If we also assume that $T\geq 0$ this now implies that $T=0$.
  \proofend

The following is the main result of this section:

  \state Theorem
  \label \absorption 
  Given any representation $\Pi$ of\/ $C^*(\B)$ on a Hilbert space $H$
there exists a representation $\Pi'$ of $\CrB$ on
$H\*\ell_2(G)$ such that the diagram 
  \def\bigbox#1{\hbox to 7em{#1}}
  $$
  \bigbox{\hfill $C^*(\B)$}
  \quad \vflx{\Pi\*\Lambda}{4} \quad
  \bigbox{$\Lin\Big(H\*\ell_2(G)\Big)$\hfill} $$$$
  \Lambda\ \raise 4pt \hbox{$\searrow$}\hskip 3em \raise 4pt
\hbox{$\nearrow$}\ \Pi'$$$$
  \CrB
  $$
  commutes.   If moreover $\Pi$ is faithful on $B_e$ then $\Pi'$ is
faithful on $\CrB$.

  \proof
  Identifying $H\*\delta_e$ and $H$ in the obvious fashion we will
think of the compression operator $F_e$ defined above as taking values
in $H$.
  Observe that for every $b_t$ in each $B_t$ we have
  $$
  F\(\Pi\*\Lambda(b_t)\) =
  F\(\pi(b_t)\*\lambda_t \) = \left\{ \matrix{
  \pi(b_e) & \hbox{ if $t=e$,} \cr
         0 & \hbox{ if $t\neq e$.}}\right.
  $$
  It follows that 
  $$
  F\circ(\Pi\*\Lambda) = \pi\circ E,
  $$
  where $E$ is the standard conditional expectation from $C^*(\B)$ to
$B_e$ (see \scite{\Amena}{2.9}).  Given $b\in C^*(\B)$ we then have by
Lemma \lcite{\FIsFaithful} that
  $$
  (\Pi\*\Lambda)(b)= 0 \sse
  F\((\Pi\*\Lambda)(b^*b)\)= 0 \sse
  \pi\(E(b^*b)\) = 0.
  \eqno{(\dagger)}
  $$
  This implies that $\Pi\*\Lambda$ vanishes on the kernel of
$\Lambda$ (see \lcite{\Recollection}) so that $\Pi\*\Lambda$  indeed factors
through $\CrB$.  If we moreover assume that $\Pi$ is faithful on
$B_e$ we have by $(\dagger)$ that the kernel of $\Lambda$ is precisely equal to the
kernel of $\Pi\*\Lambda$ and hence that $\Pi'$ is one-to-one.
  \proofend

  Our next result is an important consequence of the above result.

  \state Corollary
  \label \WrongWay
  There exists an injective *-homomorphism 
  $$\Phi : \CrB\to C^*(\B)\* \CG$$ 
  (here and elsewhere $\*$ denotes the \stress{minimal} tensor product
of C*-algebras)   such that $\Phi(b_t) = b_t\*\lambda_t$ for every $b_t$ in each $B_t$.

  \proof
  Pick a faithful representation $\Pi$ of $C^*(\B)$ on a Hilbert space
$H$.  Then the representation $\Pi'$ given by \lcite{\absorption}
provides the desired map.
  \proofend

A last technical result is in order before we embark on our main
result.

  \state Lemma
  \label \LemaDoPsi
  There exists a bounded linear map 
  $$
  \Psi : C^*(\B)\*\CG \to \CrB
  $$
  such that for every $t,s\in G$ and every $b_t \in B_t$ one has that
  $\displaystyle
  \Psi(b_t\*\lambda_s) =
  \left\{ \matrix{
  b_t & \hbox{ if $t=s$,} \cr
         0 & \hbox{ if $t\neq s$.}}\right.
  $

  \proof
  Let $\Pi$ be a faithful representation of $C^*(\B)$ on a Hilbert
space $H$ so that $\Pi\*\Lambda$ is a representation of $C^*(\B)$ on
$H\* \ell_2(G)$.  Denote by $\lambda$ the identical representation of
$\CG$ on $\ell_2(G)$ so that $(\Pi\*\Lambda)\*\lambda$ is a
representation of $C^*(\B)\*\CG$ on $H\* \ell_2(G) \* \ell_2(G)$ as
seen above.  Consider the subspace $K$ of the latter given by
  $$
  K = \bigoplus_{t\in G} H\*\delta_t\*\delta_t,
  $$
  and let $P$ be the orthogonal projection onto $K$.  Also consider the map
  $$
  \Psi: C^*(\B)\*\CG \to \Lin(K)
  $$
  given by
  $\Psi(x) = P\Big( \((\Pi\*\Lambda)\*\lambda\) (x) \Big)P$.

  There is an obvious isometric isomorphism between $K$ and
$H\*\ell_2(G)$ under which a vector of the form
$\xi\*\delta_t\*\delta_t$ is mapped to $\xi\*\delta_t$.  If we
identify $\Lin(K)$ with $\Lin(H\*\ell_2(G))$ under this map one sees that
  $$
  \Psi(b_t\*\lambda_s) =
  \left\{ \matrix{
  \Pi(b_t)\*\lambda_t & \hbox{ if $t=s$,} \cr
         0 & \hbox{ if $t\neq s$.}}\right.
  $$
  for every $t,s\in G$ and every $b_t \in B_t$.
  The range of $\Psi$ is therefore the closed linear span of
  the set of all $\Pi(b_t)\*\lambda_t$
  which happens to be also the range of $\Pi\*\Lambda$, which in turn
is isomorphic to $\CrB$ by \lcite{\absorption}.  We may then view
$\Psi$ as taking values in $\CrB$, thus providing the desired map.
  \proofend

  \section{Reduced cross-sectional algebras and exact groups}
  Let $\B$ be a Fell bundle over the discrete group $G$ and let  $\J$
be an ideal in $\B$.  As seen in \lcite{\FullExactSeq} there is a
natural exact sequence of full cross-sectional C*-algebras
  $$
  0 \to C^*(\J) \flx \iota C^*(\B) \flx \kappa C^*(\B/\J) \to 0
  \eqno{(\seqnumbering)}
  \label \FirstSequence
  $$

We would like to show that one gets out of this a (not necessarily
exact!) sequence of reduced cross-sectional C*-algebras.  For this we
need the following:

  \state Lemma  
  Given Fell bundles $\B$ and $\B'$ over the same discrete group $G$
let $\phi:C^*(\B)\to C^*(\B')$ be a *-homomorphism preserving fibers,
i.e., such that $\phi(B_t)\subseteq B'_t$ for all $t$ in $G$.  Then
there exists a *-homomorphism $\phi_r:\CrB\to \Cr{\B'}$ such that the
diagram 
  $$
  \matrix{
  C^*(\B) & \flx\phi & C^*(\B')\cr\cr
  \Lambda \flxb && \flxb \Lambda' \cr\cr
  C_r^*(\B) & \flx{\phi_r} & C_r^*(\B')
  }
  $$
  commutes, where $\Lambda$ and $\Lambda'$ are as in
\lcite{\Recollection} relative to $\B$ and $\B'$ respectively.
  
  \proof 
  Consider the diagram
  $$
  \matrix{
  C^*(\B) & \flx\phi & C^*(\B')\cr\cr
  E \flxb && \flxb E'\cr\cr
  \B_e  & \flx{\phi|_{B_e}} & B_e'
  }
  $$
  where $E$ and $E'$ are the respective conditional expectations.
This is clearly commutative and hence if $b\in C^*(\B)$ is such that
$E(b^*b)=0$ then by \lcite{\Recollection} we have that
$E'(\phi(b)^*\phi(b))=0$.  It follows that $\phi$ sends the kernel of
$\Lambda$ into the kernel of $\Lambda'$ from which the existence of
$\phi_r$ follows.
  \proofend
  
Given the exact sequence \lcite{\FirstSequence} we therefore get the (not
necessarily
exact!) sequence 
  $$ 
  0 \to \Cr\J \flx{\iota_r} \CrB \flx{\kappa_r} \Cr{\B/\J} \to 0.
  \eqno{(\seqnumbering)}
  \label \SecondSequence
  $$
  Our main result is precisely related to the situations in which one
may guarantee such a sequence to be exact:

  \state Theorem
  \label \MainTheorem
  The following conditions on a discrete group $G$ are equivalent:
  \izitem
  \zitem $G$ is an exact group.
  \zitem For every Fell bundle $\B$ over $G$ and every ideal $\J$ of
$\B$ one has that the sequence \lcite{\SecondSequence} is exact.

  \proof
  Consider the diagram 
  $$ 
  \matrix{
  0 & \to & C^*(\J)\*\CG & \vflx{\iota\*1}{3} & C^*(\B)\*\CG &
\vflx{\kappa\*1}{3} & C^*(\B/\J)\*\CG & \to & 0 \cr\cr
  &&\Phi_1 \flxc \qquad   
  &&\Phi_2 \flxc \qquad   
  &&\Phi_3 \flxc \quad
  \cr\cr
  0 & \to & \Cr\J & \vflx{\iota_r}{2.5} & \CrB & \vflx{\kappa_r}{2.5}
& \Cr{\B/\J} & \to & 0
  }
  $$
  where the vertical arrows are from \lcite{\WrongWay}.  It is a
simple matter to show that this is a commutative diagram.  

Assuming that $G$ is an exact group we have by \lcite{\FullExactSeq}
that the top row is exact and we must then prove that the same holds
for the bottom one.  The only nontrivial aspect of doing so is showing
that the kernel of $\kappa_r$ is contained in the image of $\iota_r$
which we now set out to prove.  So let $b\in \CrB$ be such that
$\kappa_r(b)=0$.  By a standard diagram chasing argument one concludes
that there exists $c\in C^*(\J)\*\CG$ such that $(\iota\*1)(c) =
\Phi_2(b)$.

We will now consider the maps $\Psi_1$ and $\Psi_2$ provided by
\lcite{\LemaDoPsi} both for $\J$ and for $\B$.  These are easily seen
to form the commutative diagram
  $$
  \matrix{
  C^*(\J)\*\CG & \vflx{\iota\*1}{3} & C^*(\B)\*\CG \cr\cr
  \Psi_1 \flxb \qquad   
  &&\Psi_2 \flxb \qquad
  \cr\cr
  \Cr\J & \vflx{\iota_r}{2.5} & \CrB 
  }
  $$
  It is also easy to see that the compositions $\Psi_i\circ\Phi_i$
give the identity maps in both cases.  We then have that
  $$
  \Phi_2 (\iota_r (\Psi_1(c))) = 
  \Phi_2 ( \Psi_2 ( \iota\*1(c))) = 
  \Phi_2 ( \Psi_2 ( \Phi_2(b))) =
  \Phi_2(b).
  $$
  Since $\Phi_2$ is injective by \lcite{\WrongWay} we conclude that
  $\iota_r (\Psi_1(c)) = b$ and hence that $b$ is in the range of
$\iota_r$.  The remaining points in the proof of the exactness of our
sequence are left to the reader.

Conversely, if (ii) is assumed, let 
  $$
  0 \to J \to B \to B/J \to 0
  $$
  be an exact sequence of C*-algebras.  If $A$ denotes any one of the
above three C*-algebras, 
consider the ``trivial Fell bundle'' $A\x G$ equipped with the
operations
  $(a,t)(b,s) = (ab,ts)$ and $(a,t)^*=(a^*,t\inv)$, 
  for $a,b\in A$ and $t,s\in G$.  It is then easy to see that $J\x G$
is an ideal in $B\x G$ and that the sequence \lcite{\SecondSequence},
which by hypothesis is exact, becomes
  $$
  0 \to J\*\CG \to B\*\CG \to B/J\*\CG \to 0.
  $$ 
  Thus $G$ is an exact group.
  \proofend

  \section{Induced and invariant ideals}
  Recall from \scite{\Amena}{3.1} that, given a discrete group $G$, a
C*-algebra $B$ is said to be \stress{graded} (over $G$) when $B$ contains a
linearly independent collection of closed subspaces $\{B_t\}_{t\in G}$
such that for each $t,s\in G$ one has
  \izitem
  \zitem $B_t^* = B_{t\inv}$
  \zitem $B_t B_s \subseteq B_{ts}$
  \zitem $B= \overline{\bigoplus_{t\in G} B_t}.$ 

  \medskip \noindent According to \scite{\Amena}{3.4} such a grading
is said to be a \stress{topological} grading if there exists a bounded
linear map $F_e: B\to B_e$ which restricts to the identity on $B_e$
and which vanishes on all other $B_t$'s.  If follows that $F_e$ is a
positive contractive conditional expectation onto $B_e$
\scite{\Amena}{3.3} and that there are contractive projections 
  $$
  F_t: B \to B_t
  $$
  such that for all finite sums $x=\sum_{t\in G} b_t$, with $b_t \in
B_t$, one has $F_t(x) = b_t$ \scite{\Amena}{3.5}.  Throughout this
section we shall fix a topologically graded C*-algebra $B=
\overline{\bigoplus_{t\in G} B_t}.$

Given an ideal $\Ideal$ in $B$ one may consider the following sets:
  \iItem
  \Item{} $\Ideal_1=\<\Ideal \cap B_e \>$, i.e.~the ideal generated by
$\Ideal \cap B_e$,
  \Item{}$\Ideal_2=\{ b\in B: F_e(b^* b) \in \Ideal  \}$, and
  \Item{} $\Ideal_3=\{ b\in B: F_t(b) \in \Ideal, \hbox{ for all }
t\in G\}$.

  \medskip\noindent One always has $\Ideal_1 \subseteq \Ideal_2 =
\Ideal_3$ \scite{\Amena}{3.9} and whenever the underlying Fell bundle
$\B = \{B_t\}_{t\in G}$ has the \stress{approximation property}
defined in \scite{\Amena}{4.4}, e.g.~if $G$ is amenable
\scite{\Amena}{4.7}, one also has that $\Ideal_1=\Ideal_2$
\scite{\Amena}{4.10}.
  This should be compared to similar results by Str\v{a}til\v{a} and
Voiculescu on AF-algebras \cite{\SV} and by Nica on quasi-lattice
ordered groups \cite{\Nica}.
  It is our goal in this section to exhibit still another situation in
which one may guarantee that $\Ideal_1=\Ideal_2$.

  \state Theorem
  \label \MainApplication
  Let $G$ be a discrete group and let 
  $B = \overline{\bigoplus_{t\in G} B_t}$ be a topologically
$G$-graded C*-algebra.  Suppose that
  $G$ is exact and that
  the standard conditional expectation $F:B\to B_e$ is
faithful.  Then for every ideal $\Ideal$ of $B$ one has that the
sets $\Ideal_1$ and $\Ideal_2$ defined above coincide.

  \proof 
  Denote by $\B = \{B_t\}_{t\in G}$ the underlying Fell bundle and
note that by \scite{\Amena}{3.7} $B$ is isomorphic to $\CrB$.  For
each $t\in G$\/ let $J_t=\Ideal\cap B_t$ so that $\J:=\{J_t\}_{t\in G}$ is an
ideal in $\B$ according to \lcite{\DefineIdeal}.  By
\lcite{\MainTheorem} the sequence
  $$
  0 \to \Cr{\J} \flx{\iota_r} B \flx{\kappa_r} \Cr{\B/\J} \to 0
  $$
  is exact.  We now claim that 
  $\Ideal_1 = \iota_r\(\Cr{\J}\)$ and $\Ideal_2 = \Ker(\kappa_r)$ from
where the proof will be concluded.

  For each $t$ in $G$ we have by \scite {\BMS}{1.7} that 
  $$
  J_t = J_t J_t^* J_t \subseteq J_e J_t \subseteq \Ideal_1
  $$
  so that $\iota_r\(\Cr{\J}\) =\Ideal_1$.  On the other hand, denoting
by $E$ the standard conditional expectation of $\Cr{\B/\J}$ we have
for any $b\in B$ that
  $$
  \kappa_r(b) = 0 \sse
  E(\kappa_r(b^*b)) = 0 \sse
  \kappa_r(F_e(b^*b)) = 0 \sse
  F_e(b^*b) \in \Ideal.
  \proofend
  $$

According to \scite{\Amena}{3.10} an ideal $\Ideal$ in $B$ is said to
be an \stress{induced} ideal when $\Ideal = \<\Ideal \cap B_e\>$.  Let
us also make the following:

  \definition
  An ideal $\Ideal$ in $B$ is said to be \stress{invariant} if
$F_t(\Ideal)\subseteq \Ideal$ for all $t$ in $G$.

  \state Corollary
  Assume that $G$ is exact and that the standard conditional
expectation $F_e$ is faithful on $B$.  Then any invariant ideal is
induced.

  \proof
  Let $\Ideal$ be an invariant ideal in $B$.  Then
  $$
  \Ideal \subseteq 
  \big\{ b\in B: F_t(b) \in \Ideal, \hbox{ for all } t\in G \big\} = 
  \Ideal_3 = \Ideal_1 = 
  \<\Ideal \cap B_e \> \subseteq
  \Ideal,
  $$
  so equality holds throughout.
  \proofend

  \endgroup
  \end